\documentclass[11pt,letterpaper]{article}

\usepackage[margin=1in]{geometry}
\usepackage{amsmath,amssymb,amsthm,mathtools}
\usepackage{algorithm}
\usepackage{float}
\usepackage{algpseudocode}
\usepackage{microtype}
\usepackage{enumitem}
\usepackage{booktabs}
\usepackage[hidelinks]{hyperref}
\usepackage[T1]{fontenc}
\usepackage{lmodern}

\newtheorem{theorem}{Theorem}[section]
\newtheorem{lemma}[theorem]{Lemma}
\newtheorem{proposition}[theorem]{Proposition}

\theoremstyle{definition}
\newtheorem{definition}[theorem]{Definition}
\theoremstyle{remark}
\newtheorem{remark}[theorem]{Remark}

\numberwithin{equation}{section}

\newcommand{\Z}{\mathbb Z}
\newcommand{\ph}{\varphi}
\newcommand{\floor}[1]{\left\lfloor #1\right\rfloor}

\title{\textbf{On the Frobenius Coin Problem in Three Variables:\\
A New Algorithm and Some Analysis}}
\author{
Negin Bagherpour\thanks{Department of Engineering Sciences, University of Tehran,
Tehran, Iran. E-mail: \texttt{negin.bagherpour@ut.ac.ir}}
\and
Amir Jafari\thanks{Department of Mathematical Sciences, Sharif University of Technology,
Tehran, Iran.}
\and
Amin Najafi\footnotemark[2]
}
\date{}

\begin{document}
\maketitle

\begin{abstract}
The three-variable Frobenius coin problem asks for the largest integer that
cannot be represented as a non-negative integer combination of three positive
integers $a_1<a_2<a_3$ with $\gcd(a_1,a_2,a_3)=1$. We introduce an
algorithmic framework based on two synchronized Euclidean-division sequences.
The method has arithmetic complexity $O(\log a_1)$, matching the best known
logarithmic order for the three-variable problem. We prove that this
dependence is asymptotically sharp for the proposed algorithm by constructing
an explicit infinite Fibonacci-type family for which Step~2 requires
$\frac12\log_\ph(a_1)+O(1)$ iterations. Hence the worst-case arithmetic
complexity of the algorithm is $\Theta(\log a_1)$. We also give a
constant-overhead reduction from gcd computation to two evaluations of the
three-variable Frobenius function. In addition, the synchronized framework
yields an upper bound for the Frobenius number, recovers several known
special cases, and provides another proof of the Erd\H{o}s--Graham bound.
The main advantage of the method is its elementary Euclidean structure,
which supports both computation and theoretical analysis.
\end{abstract}



\section{Introduction}

According to Brauer~\cite{Brauer1942}, Frobenius mentioned the following
problem in his lectures in the early twentieth century. For relatively prime
positive integers
\[
a_1<a_2<\cdots<a_n,
\]
find the largest integer $G(a_1,\ldots,a_n)$ that cannot be represented as
\[
a_1x_1+\cdots+a_nx_n
\]
with positive integer coefficients $x_1,\ldots,x_n$. If
$g(a_1,\ldots,a_n)$ denotes the largest integer that cannot be represented
with non-negative integer coefficients, then
\begin{equation}\label{eq:G-g-relation}
g(a_1,\ldots,a_n)
=
G(a_1,\ldots,a_n)-a_1-\cdots-a_n.
\end{equation}
The usual ``coin problem'' interpretation asks for the largest monetary
amount that cannot be formed from an unlimited supply of coins with
denominations $a_1,\ldots,a_n$.

The quantities $g$ and $G$ will both be called Frobenius numbers, with the
convention always made clear from the context. For $n=2$, Sylvester proved
\[
G(a_1,a_2)=a_1a_2,
\qquad
g(a_1,a_2)=a_1a_2-a_1-a_2.
\]
For $n>2$ no comparably simple universal closed formula is known; see, for
example, Curtis~\cite{Curtis1990}. Tripathi~\cite{Tripathi2017} gave
algorithmic formulae for the three-variable case. To make the formulation of Tripathi \cite{Tripathi2017} practically applicable—which inherently involves non-trivial subproblems—Suhajda and
Thillaisundaram~\cite{Suhajda2026} recently introduced a case-by-case classification. For completeness, we begin with the standard fact that the Frobenius number
is well defined.

\begin{lemma}\label{lem:well-defined}
Let $a_1,\ldots,a_n$ be positive integers with
$\gcd(a_1,\ldots,a_n)=1$. There is an integer $m_0$ such that every
$m\ge m_0$ can be written as
\[
m=a_1x_1+\cdots+a_nx_n,
\qquad x_1,\ldots,x_n\in\Z_{\ge0}.
\]
\end{lemma}

\begin{proof}
By B\'ezout's theorem, choose integers $x_1,\ldots,x_n$ such that
\[
a_1x_1+\cdots+a_nx_n=1.
\]
After reordering if necessary, assume
$x_1,\ldots,x_i\ge0$ and $x_{i+1},\ldots,x_n<0$. Let
\[
P=a_1x_1+\cdots+a_ix_i,
\qquad
N=-a_{i+1}x_{i+1}-\cdots-a_nx_n.
\]
Then $P-N=1$. We claim that $m_0=(a_1-1)N$ works. For $m\ge m_0$, write
\[
m-(a_1-1)N=a_1q+r,
\qquad 0\le r<a_1.
\]
Then
\[
m=(a_1-1-r)N+a_1q+rP,
\]
which is a non-negative integer combination of the $a_i$.
\end{proof}

We shall also use the following elementary uniqueness statement.

\begin{lemma}\label{lem:two-var-representation}
If $\gcd(a_1,a_2)=1$, then every positive integer $n$ can be written uniquely
as
\[
n=a_1x_1+a_2x_2,
\qquad x_1,x_2\in\Z,
\qquad 0<x_1\le a_2.
\]
\end{lemma}

\begin{proof}
By B\'ezout's theorem, choose integers $y_1,y_2$ such that
$a_1y_1+a_2y_2=n$. Write
$y_1=qa_2+x_1$ with $0<x_1\le a_2$. Then
\[
n=a_1x_1+a_2(a_1q+y_2),
\]
which proves existence.

For uniqueness, suppose
\[
a_1x_1+a_2x_2=a_1x'_1+a_2x'_2
\]
with $0<x_1,x'_1\le a_2$. Then
\[
a_1(x_1-x'_1)=a_2(x'_2-x_2).
\]
Since $\gcd(a_1,a_2)=1$, the integer $x_1-x'_1$ is divisible by $a_2$.
But $|x_1-x'_1|<a_2$, hence $x_1=x'_1$ and then $x_2=x'_2$.
\end{proof}

\begin{remark}\label{rem:sylvester}
Sylvester's formula follows quickly from
Lemma~\ref{lem:two-var-representation}. In particular,
$G(a_1,a_2)=a_1a_2$, or equivalently
$g(a_1,a_2)=a_1a_2-a_1-a_2$.
\end{remark}

\begin{remark}\label{rem:euclid-two-var}
The representation in Lemma~\ref{lem:two-var-representation} can be
constructed by successive divisions, as in the Euclidean algorithm. The
classical Fibonacci analysis of the Euclidean algorithm gives a logarithmic
bound in the smaller input; see Knuth~\cite{Knuth1998}.
\end{remark}

The main object of this paper is an algorithm for
$G(a_1,a_2,a_3)$ when
\[
0<a_1<a_2<a_3,
\qquad
\gcd(a_1,a_2,a_3)=1.
\]
Without loss of generality, the algorithmic core may be reduced to the case
$\gcd(a_1,a_2)=1$ by the following result of
Johnson~\cite{Johnson1960}.

\begin{lemma}[Johnson]\label{lem:johnson}
Suppose $\gcd(a_1,a_2,a_3)=1$ and let $d=\gcd(a_1,a_2)$. Then
\[
G(a_1,a_2,a_3)
=
d\,G\!\left(\frac{a_1}{d},\frac{a_2}{d},a_3\right).
\]
\end{lemma}

Thus a preliminary gcd computation reduces the problem to
$\gcd(a_1,a_2)=1$.

\paragraph{Arithmetic complexity model.}
Throughout this paper, complexity refers to \emph{arithmetic complexity},
not bit complexity. We use a unit-cost integer arithmetic model in which
additions, subtractions, multiplications, comparisons, and Euclidean
divisions are counted as elementary arithmetic operations independently of
the bit length of their operands. In particular, one division
\[
x=qy+r,\qquad 0\le r<y,
\]
counts as one division step.

The three-variable problem has a substantial algorithmic literature.
Selmer and Beyer~\cite{SelmerBeyer1978} developed a continued-fraction
approach, and R{\o}dseth~\cite{Rodseth1979} gave a related refinement.
Davison~\cite{Davison1994} subsequently presented a polynomial-time
algorithm with similarities to Euclidean/continued-fraction procedures.
Importantly for the complexity comparison made here,
Greenberg~\cite{Greenberg1988} already obtained logarithmic behavior in the
smallest coefficient. Consequently, the $O(\log a_1)$ estimate proved below
should \emph{not} be interpreted as a new asymptotic order better than every
previous method. The distinguishing features of the present algorithm are
its particularly elementary synchronized Euclidean structure and the way
this structure can be reused in the subsequent analysis.

Tripathi's formulae~\cite{Tripathi2017} also deserve a computational
comment. Several of the formulae contain an auxiliary step in which, for
rational numbers $0<\alpha<\beta$, one has to determine the largest integer
$m$ such that
\[
\floor{m\alpha}=\floor{m\beta}.
\]
The values of $m$ satisfying this equality need not be consecutive. For
example, for
\[
a_1=468342493,\qquad
a_2=472518070,\qquad
a_3=472714471,
\]
our synchronized algorithm produces
\[
\begin{array}{c|rrrrrr}
 & -1&0&1&2&3&4\\ \hline
 b_i
 &a_2&235128713&2260644&21737&21733&4\\
 b'_i
 &a_1&233050914&2240665&21754&3&21751
\end{array}
\]
and returns
\[
G(a_1,a_2,a_3)
=
b_3a_1+b'_4a_2-\min\{b_4a_1,b'_3a_2\}
=
20454810386729.
\]
At one of the corresponding stages in the formulae of
Tripathi~\cite{Tripathi2017}, one encounters
\[
\alpha=
\frac{10906201339369943}{105875470131655},
\qquad
\beta=\frac{232868069}{2260644}.
\]
For $m\ge10400$, the relevant values satisfying
$\floor{m\alpha}=\floor{m\beta}$ include
\[
10400,10401,10402,10403,10504,10505,10506,10608,10609,10712,
\]
so the largest one is $10712$. This illustrates that an explicit formula
can still contain a nontrivial auxiliary search unless that step is supplied
with its own effective implementation.

We now present the algorithm.

\section{The algorithm and its correctness}

Assume throughout this section that
\[
0<a_1<a_2<a_3,
\qquad
\gcd(a_1,a_2)=1.
\]

\begin{algorithm}[H]
\caption{Three-variable Frobenius coin problem solver}
\label{alg:main}
\begin{algorithmic}[1]
\Require Integers $0<a_1<a_2<a_3$ with $\gcd(a_1,a_2)=1$.
\Ensure $G(a_1,a_2,a_3)$.
\State \textbf{Step 0.} Set $b'_{-1}=a_1$ and $b_{-1}=a_2$.
\State \textbf{Step 1.} Find $0<b'_0\le a_1$ and $b_0$ such that
\[
a_3=a_2b'_0-a_1b_0.
\]
\If{$b_0\le0$}
\State \Return $G(a_1,a_2,a_3)=a_1a_2+a_3$.
\EndIf
\State \textbf{Step 2.}
\For{$i=0,1,\ldots$}
\If{$b_i=0$ or $b'_i=0$}
\State \textbf{break}
\EndIf
\State Compute
\[
b_{i-1}=b_iq_i+r_i,
\qquad
0\le r_i<b_i,
\]
and
\[
b'_{i-1}=b'_iq'_i+r'_i,
\qquad
0\le r'_i<b'_i.
\]
\State Set $c_i=\min(q_i,q'_i)$ and compute
\[
b_{i+1}=b_{i-1}-c_ib_i,
\qquad
b'_{i+1}=b'_{i-1}-c_ib'_i.
\]
\If{$q_i\ne q'_i$}
\State \textbf{break}
\EndIf
\EndFor
\State \textbf{Step 3.} Let $k$ be the final computed index.
\If{$k$ is odd}
\State \Return
\[
b_ka_1+b'_{k-1}a_2-\min\{b_{k-1}a_1,b'_ka_2\}.
\]
\Else
\State \Return
\[
b_{k-1}a_1+b'_ka_2-\min\{b_ka_1,b'_{k-1}a_2\}.
\]
\EndIf
\end{algorithmic}
\end{algorithm}

A few immediate observations are useful.

First, the numbers $b_0,b'_0$ in Step~1 are uniquely determined by
Lemma~\ref{lem:two-var-representation}. If $b_0\le0$, then $a_3$ is a
non-negative integer combination of $a_1$ and $a_2$, hence it is redundant:
\[
g(a_1,a_2,a_3)=g(a_1,a_2),
\]
and therefore
\[
G(a_1,a_2,a_3)=a_1a_2+a_3.
\]
Since $a_3>0$, one also has $b_0\le a_2$.

Second, before termination the quantities $b_{i+1}$ and $b'_{i+1}$ are the
actual Euclidean remainders whenever $q_i=q'_i$. Thus
$b_{i+1}<b_i$ and $b'_{i+1}<b'_i$, so the loop terminates. A refined
iteration bound is proved in Section~\ref{sec:complexity}.

Finally, the same Euclidean pass used to compute $\gcd(a_1,a_2)$ can be
reused to obtain the initial representation in Step~1. Thus all preparatory
work is logarithmic in $a_1$ in the arithmetic model.

We now prove correctness. Suppose the algorithm terminates after the final
pair $b_k,b'_k$ has been computed.

\begin{lemma}\label{lem:p-identity}
For $-1\le j\le k$,
\begin{equation}\label{eq:p-identity}
b'_ja_2=b_ja_1+(-1)^jp_ja_3,
\end{equation}
where
\[
p_{-1}=0,\qquad p_0=1,
\]
and
\begin{equation}\label{eq:p-recurrence}
p_{j+1}=\min(q_j,q'_j)p_j+p_{j-1}
\end{equation}
whenever the next term is defined. Moreover, $p_j\ge F_j$, where
$F_0=0,F_1=1$. If all common quotients before the $j$th term are equal to
$1$, then the corresponding $p_i$ are Fibonacci numbers; if at some stage
the common quotient is at least $2$, the Fibonacci lower bound becomes
strict at the next relevant index.
\end{lemma}

\begin{proof}
For $j=0$, identity \eqref{eq:p-identity} is exactly the defining relation
\[
a_3=a_2b'_0-a_1b_0,
\]
and therefore $p_0=1$. The identity also holds for $j=-1$ with
$p_{-1}=0$. Assume it holds through index $i$. Subtracting
$\min(q_i,q'_i)$ times the identity at index $i$ from the identity at index
$i-1$ gives
\[
b'_{i+1}a_2
=
b_{i+1}a_1
+
(-1)^{i+1}
\bigl(p_{i-1}+\min(q_i,q'_i)p_i\bigr)a_3,
\]
which proves both \eqref{eq:p-identity} and
\eqref{eq:p-recurrence}. Since every common quotient used in the recurrence
is at least $1$, the Fibonacci lower bound follows by induction. In
particular, if a common quotient is at least $2$, then
\[
p_{j+1}\ge2F_j+F_{j-1}=F_{j+2}.
\]
\end{proof}

\begin{lemma}\label{lem:determinants}
For every $0\le j\le k$,
\begin{align}
a_1&=p_jb'_{j-1}+p_{j-1}b'_j,\label{eq:a1det}\\
a_2&=p_jb_{j-1}+p_{j-1}b_j,\label{eq:a2det}\\
(-1)^ja_3&=b_{j-1}b'_j-b_jb'_{j-1}.\label{eq:a3det}
\end{align}
\end{lemma}

\begin{proof}
For $j=0$ the identities follow from the definitions. The general case is
obtained by induction from the recurrences defining
$b_{j+1},b'_{j+1}$ and $p_{j+1}$.
\end{proof}

\begin{lemma}\label{lem:terminal-order}
Assume the algorithm stops because the two quotients differ and
$b_kb'_k\ne0$. If $k$ is odd, then
\[
q'_{k-1}<q_{k-1},
\qquad
b'_k<b'_{k-1},
\qquad
b_k>b_{k-1}.
\]
If $k$ is even, all three inequalities are reversed.
\end{lemma}

\begin{proof}
Suppose $k$ is odd. By \eqref{eq:a3det} at index $k-1$,
\[
b'_{k-1}b_{k-2}-b'_{k-2}b_{k-1}=a_3>0,
\]
and hence
\[
\frac{b_{k-2}}{b_{k-1}}
>
\frac{b'_{k-2}}{b'_{k-1}}.
\]
Since the first mismatch occurs at the division producing $b_k,b'_k$, we
obtain
$q_{k-1}>q'_{k-1}$. Thus
\[
b'_k=b'_{k-2}-q'_{k-1}b'_{k-1}=r'_{k-1}<b'_{k-1},
\]
whereas
\[
\begin{aligned}
b_k
&=b_{k-2}-q'_{k-1}b_{k-1}\\
&\ge b_{k-2}-(q_{k-1}-1)b_{k-1}\\
&=r_{k-1}+b_{k-1}
>b_{k-1}.
\end{aligned}
\]
The case of even $k$ is symmetric.
\end{proof}

\begin{definition}\label{def:L}
For relatively prime positive integers $a_1,a_2,a_3$, let
$L(a_1\mid a_2,a_3)$ denote the smallest positive integer $m$ such that
\[
ma_1=a_2x+a_3y
\]
for some $x,y\in\Z_{\ge0}$.
\end{definition}

\begin{lemma}\label{lem:L-values}
If $k$ is odd, then
\[
b_k=L(a_1\mid a_2,a_3),
\qquad
b'_{k-1}=L(a_2\mid a_1,a_3).
\]
If $k$ is even, then
\[
b_{k-1}=L(a_1\mid a_2,a_3),
\qquad
b'_k=L(a_2\mid a_1,a_3).
\]
\end{lemma}

\begin{proof}
Assume $k$ is odd. By Lemma~\ref{lem:p-identity},
\begin{align}
b_ka_1&=b'_ka_2+p_ka_3,\label{eq:Lodd1}\\
b'_{k-1}a_2&=b_{k-1}a_1+p_{k-1}a_3.\label{eq:Lodd2}
\end{align}
Hence, with
\[
L_1=L(a_1\mid a_2,a_3),
\qquad
L_2=L(a_2\mid a_1,a_3),
\]
we have $L_1\le b_k$ and $L_2\le b'_{k-1}$.

Suppose first that $L_2<b'_{k-1}$ and
\[
L_2a_2=xa_1+ya_3,
\qquad x,y\ge0.
\]
Any integral relation among $a_1,a_2,a_3$ is generated by two consecutive
relations in Lemma~\ref{lem:p-identity}. Indeed, every later relation is
obtained from the previous two, so it is enough to verify this for the first
two:
\[
b_{-1}a_1=b'_{-1}a_2,
\qquad
b_0a_1=b'_0a_2+a_3.
\]
If $a_1x_1+a_2x_2+a_3x_3=0$, subtracting $x_3$ times the second relation
eliminates $a_3$. Since $\gcd(a_1,a_2)=1$, the remaining coefficients are
multiples of $a_2$ and $a_1$, and the first relation eliminates them.

Consequently there exist $\alpha,\beta\in\Z$ such that
\begin{align}
-\alpha b_k+\beta b_{k-1}&=x\ge0,\label{eq:Lalpha1}\\
0<-\alpha b'_k+\beta b'_{k-1}
&=L_2<b'_{k-1},\label{eq:Lalpha2}\\
\alpha p_k+\beta p_{k-1}&=y\ge0.\label{eq:Lalpha3}
\end{align}
If the algorithm stops because one of $b_k,b'_k$ is zero, then for odd
$k$ relation \eqref{eq:Lodd1} forces $b'_k=0$, and
\eqref{eq:Lalpha2} is impossible. Hence assume $b'_k\ne0$. By
Lemma~\ref{lem:terminal-order},
\[
b_k>b_{k-1},
\qquad
b'_k<b'_{k-1}.
\]
If $\alpha<0$, then \eqref{eq:Lalpha3} gives $\beta\ge0$, contradicting
\eqref{eq:Lalpha2}. Therefore $\alpha,\beta>0$.
Equation~\eqref{eq:Lalpha1} and $b_k>b_{k-1}$ imply
$\alpha<\beta$, while \eqref{eq:Lalpha2} and
$b'_k<b'_{k-1}$ imply $\alpha\ge\beta$, a contradiction. Thus
$L_2=b'_{k-1}$.

Now suppose $L_1<b_k$ and write
\[
L_1a_1=x'a_2+y'a_3,
\qquad x',y'\ge0.
\]
As above, there exist $\alpha,\beta\in\Z$ such that
\begin{align}
0<\alpha b_k-\beta b_{k-1}
&=L_1<b_k,\label{eq:Lbeta1}\\
0<\alpha b'_k-\beta b'_{k-1}
&=x'\ge0,\label{eq:Lbeta2}\\
\alpha p_k+\beta p_{k-1}
&=y'\ge0.\label{eq:Lbeta3}
\end{align}
If $b'_k=0$, then \eqref{eq:Lbeta2} gives $\beta<0$, while
\eqref{eq:Lbeta3} gives $\alpha>0$, contradicting
\eqref{eq:Lbeta1}. If $b'_k\ne0$, Lemma~\ref{lem:terminal-order} gives
$b_k>b_{k-1}$ and $b'_k<b'_{k-1}$. Then
\eqref{eq:Lbeta1} forces $\alpha\le\beta$, whereas
\eqref{eq:Lbeta2} forces $\alpha>\beta$, again a contradiction.
Therefore $L_1=b_k$.

The even case is obtained by the same argument with the roles reversed.
\end{proof}

\begin{lemma}\label{lem:max-a3}
Suppose the algorithm terminates after the final index $k$.

If $k$ is odd, then the representation
\[
b_ka_1=b'_ka_2+p_ka_3
\]
has maximal coefficient of $a_3$ among all non-negative representations of
$b_ka_1$ by $a_2$ and $a_3$. If $b'_k\ne0$, the same is true of
\[
b'_{k-1}a_2=b_{k-1}a_1+p_{k-1}a_3.
\]
If $k$ is even, the analogous assertions hold for
\[
b'_ka_2=b_ka_1+p_ka_3
\]
and, when $b_k\ne0$,
\[
b_{k-1}a_1=b'_{k-1}a_2+p_{k-1}a_3.
\]
\end{lemma}

\begin{proof}
We prove the odd case. By the proof of
Lemma~\ref{lem:L-values}, every relation among $a_1,a_2,a_3$ is generated
by the two terminal relations. To improve the coefficient of $a_3$ in
\[
b_ka_1=b'_ka_2+p_ka_3,
\]
we would need integers $\alpha,\beta$ satisfying
\[
\alpha b_k-\beta b_{k-1}=b_k,
\]
\[
\alpha b'_k-\beta b'_{k-1}\ge0,
\]
and
\[
\alpha p_k+\beta p_{k-1}>p_k.
\]
If $b'_k=0$, the second inequality forces $\beta<0$ and hence
$\alpha>0$, which is incompatible with the first equation. If $b'_k\ne0$,
Lemma~\ref{lem:terminal-order} gives
$b_k>b_{k-1}$ and $b'_k<b'_{k-1}$. The first relation forces
$\alpha\le\beta$, while the second forces $\alpha>\beta$, a contradiction.

For the second terminal representation, increasing the coefficient of
$a_3$ would require
\[
-\alpha b'_k+\beta b'_{k-1}=b'_{k-1},
\]
\[
-\alpha b_k+\beta b_{k-1}\ge0,
\]
and
\[
\alpha p_k+\beta p_{k-1}>p_{k-1}.
\]
When $b'_k\ne0$, if $\alpha<0$ then the first equation gives
$\beta>0$, contradicting that equation. Hence $\alpha,\beta\ge0$.
The second inequality and $b_k>b_{k-1}$ give $\alpha<\beta$,
whereas the first equation and $b'_k<b'_{k-1}$ give
$\alpha\ge\beta$. This contradiction proves the claim. When $b'_k=0$,
representations with a larger $a_3$ coefficient may exist, but that case is
handled separately in the terminal formula below.

The even case is symmetric.
\end{proof}

The next lemma is a form of Johnson's characterization adapted to the
present assumptions.

\begin{lemma}\label{lem:johnson-characterization}
Let
\[
0<a_1<a_2<a_3,
\qquad
\gcd(a_1,a_2)=1,
\]
and assume that $a_3$ is not a non-negative integer combination of
$a_1,a_2$. Put
\[
L_1=L(a_1\mid a_2,a_3),
\qquad
L_2=L(a_2\mid a_1,a_3),
\]
and choose
\[
L_1a_1=xa_2+ya_3,
\qquad
L_2a_2=x'a_1+y'a_3
\]
with non-negative coefficients, taking $y,y'$ maximal when the
representations are not unique. Then
\begin{align}
G(a_1,a_2,a_3)
&=
\max\{L_1a_1+y'a_3,\;L_2a_2+ya_3\}\label{eq:johnson-max}\\
&=
L_1a_1+L_2a_2-\min\{x'a_1,xa_2\}.\label{eq:johnson-min}
\end{align}
\end{lemma}

\begin{proof}
Since $\gcd(a_1,a_2)=1$, Lemma~\ref{lem:two-var-representation} gives
\[
a_3=a_1s_1+a_2s_2,
\qquad
0<s_1\le a_2.
\]
Because $a_3$ is not a non-negative combination of $a_1,a_2$, we have
$0<s_1<a_2$ and $s_2<0$. Thus $s_1a_1$ is a positive combination of
$a_2$ and $a_3$, so $L_1\le s_1<a_2$. Similarly $L_2<a_1$. In
particular, the maximal coefficients $y,y'$ above are positive.

Now $G=G(a_1,a_2,a_3)$ is not a positive integer combination of all three
generators, whereas $G+a_i$ is such a combination for each
$i=1,2,3$. Under the present assumptions, this is equivalent to $G$ having
a positive representation by each pair of generators. We show that only the
two values in \eqref{eq:johnson-max} can have these properties.

Suppose
\[
N=za_2+z'a_3
 =ta_1+t'a_3
 =wa_1+w'a_2
\]
with all six coefficients positive, and choose the first two
representations so that $z'$ and $t'$ are maximal. Assume that $N$ itself
is not a positive combination of $a_1,a_2,a_3$.

First, $w\le L_1$ and $w'\le L_2$. For example, if $w>L_1$, then
\[
\begin{aligned}
N
&=(w-L_1)a_1+L_1a_1+w'a_2\\
&=(w-L_1)a_1+(x+w')a_2+ya_3,
\end{aligned}
\]
a contradiction. The proof of $w'\le L_2$ is the same.

Similarly $z\le L_2$ and $t\le L_1$. If, for instance, $z>L_2$ and
$x'>0$, substitute
$L_2a_2=x'a_1+y'a_3$ into the representation of $N$ to obtain a positive
three-generator representation. If $x'=0$, replacing $z$ by
$z-L_2$ increases $z'$, contradicting maximality. The argument for
$t\le L_1$ is identical.

We compare $z'$ and $t'$. If $z'=t'$, then
$za_2=ta_1$, so $a_1\mid z$ and therefore
$z\ge a_1$, contradicting $z\le L_2<a_1$.

If $z'>t'$, then
\[
ta_1=za_2+(z'-t')a_3,
\]
so $t\ge L_1$. Hence $t=L_1$. It follows that
\[
w'a_2=(L_1-w)a_1+t'a_3,
\]
so $w'\ge L_2$, and therefore $w'=L_2$. By maximality of $t'$ and $y'$,
we have $t'=y'$, whence
\[
N=L_1a_1+y'a_3.
\]

If $z'<t'$, then
\[
za_2=ta_1+(t'-z')a_3,
\]
and therefore $z\ge L_2$. Thus $z=L_2$. Consequently
\[
wa_1=(L_2-w')a_2+z'a_3,
\]
so $w\ge L_1$, and hence $w=L_1$. By maximality of $z'$ and $y$,
we have $z'=y$, giving
\[
N=L_2a_2+ya_3.
\]
This proves \eqref{eq:johnson-max}. Formula
\eqref{eq:johnson-min} follows immediately by substituting the defining
representations of $L_1a_1$ and $L_2a_2$.
\end{proof}

\begin{theorem}\label{thm:correctness}
If Algorithm~\ref{alg:main} stops after the final index $k$, then
\[
G(a_1,a_2,a_3)
=
\begin{cases}
b_ka_1+b'_{k-1}a_2-\min\{b_{k-1}a_1,b'_ka_2\},
& k\ \text{odd},\\[2mm]
b_{k-1}a_1+b'_ka_2-\min\{b_ka_1,b'_{k-1}a_2\},
& k\ \text{even}.
\end{cases}
\]
Hence Algorithm~\ref{alg:main} is correct.
\end{theorem}

\begin{proof}
If $b_kb'_k\ne0$, the result follows directly from
Lemmas~\ref{lem:L-values}, \ref{lem:max-a3}, and
\ref{lem:johnson-characterization}.

It remains to explain the zero-remainder case. Suppose $k$ is odd and
$b'_k=0$. Then Lemma~\ref{lem:L-values} gives
$L_1=b_k$, and
\[
b_ka_1=p_ka_3
\]
has maximal possible coefficient of $a_3$. In the notation of
Lemma~\ref{lem:johnson-characterization}, the coefficient $x$ is zero.
Therefore the choice of representation for $L_2a_2$ does not affect the
maximum, and
\[
G(a_1,a_2,a_3)
=
L_1a_1+L_2a_2
=
b_ka_1+b'_{k-1}a_2.
\]
The even case is analogous.
\end{proof}

\section{Complexity and auxiliary lemmas}
\label{sec:complexity}

We next sharpen the Euclidean termination estimate.

\begin{lemma}\label{lem:complexity-upper}
If Step~2 of Algorithm~\ref{alg:main} terminates after $k$ iterations, then
\[
k<\frac12\log_\ph(a_1)+2,
\qquad
\ph=\frac{1+\sqrt5}{2}.
\]
\end{lemma}

\begin{proof}
We prove the odd case; the even case is analogous. With the notation of
Lemmas~\ref{lem:p-identity} and \ref{lem:determinants}, we first show
\[
p_{k-1}<\sqrt{a_1}.
\]
By \eqref{eq:a1det},
\[
a_1=p_{k-1}b'_k+p_kb'_{k-1},
\]
and Lemma~\ref{lem:p-identity} gives
\[
b'_{k-1}a_2=b_{k-1}a_1+p_{k-1}a_3,
\]
hence
\[
b'_{k-1}>
\frac{p_{k-1}a_3}{a_2}.
\]
Therefore
\[
a_1
>
p_{k-1}b'_k
+
p_k\frac{p_{k-1}a_3}{a_2}.
\]
Since $a_3>a_2$ and $p_k>p_{k-1}$,
\[
a_1>p_kp_{k-1}>p_{k-1}^2.
\]
Thus $\sqrt{a_1}>p_{k-1}\ge F_k$. For $k\ge2$, the standard elementary
Fibonacci estimate
\[
F_k\ge \ph^{\,k-2}
\]
follows by induction from $F_{j+1}=F_j+F_{j-1}$ and
$\ph^2=\ph+1$. Hence
\[
\sqrt{a_1}>\ph^{\,k-2},
\]
and therefore
\[
k<\frac12\log_\ph(a_1)+2.
\]
The cases $k=0,1$ are immediate.
\end{proof}

\begin{remark}\label{rem:fibonacci}
For every relevant index $i$, the recurrence
\eqref{eq:p-recurrence} implies a Fibonacci lower bound on $p_i$. The slowest
growth occurs when all common quotients are $1$. Theorem~\ref{thm:sharpness}
below shows that this Fibonacci regime is realized by an infinite family of
valid Frobenius instances and therefore gives the true worst-case order of
Algorithm~\ref{alg:main}.
\end{remark}

\begin{remark}\label{rem:quotient-two}
If the common quotient is $2$ over a consecutive range, the $p_i$ grow
strictly faster than in the all-$1$ Fibonacci regime. More generally, any
common quotient larger than $1$ forces a strict jump above the corresponding
Fibonacci lower bound.
\end{remark}

\begin{theorem}[Sharpness of the logarithmic bound]
\label{thm:sharpness}
The estimate in Lemma~\ref{lem:complexity-upper} is asymptotically sharp.
There exists an infinite family of admissible triples for which Step~2 of
Algorithm~\ref{alg:main} requires
\[
\frac12\log_\ph(a_1)+O(1)
\]
iterations. Consequently, the worst-case arithmetic complexity of
Algorithm~\ref{alg:main} is
\[
\Theta(\log a_1).
\]
\end{theorem}

\begin{proof}
Let $F_0=0,F_1=1$. For any integer $m\ge1$, put
\[
K=6m,
\qquad
M=F_{K+3}+1,
\]
and define
\begin{equation}\label{eq:sharp-family}
a_1=MF_{K+2}+F_{K+1},
\qquad
a_2=MF_{K+3}+F_{K+1},
\qquad
a_3=M^2.
\end{equation}
First,
\[
a_2-a_1=MF_{K+1}>0,
\]
and
\[
a_3-a_2=M-F_{K+1}=F_{K+2}+1>0,
\]
so $0<a_1<a_2<a_3$.

We also have $\gcd(a_1,a_2)=1$. Indeed,
\[
\gcd(M,F_{K+1})
=
\gcd(F_{K+3}+1,F_{K+1})
=
\gcd(F_K+1,F_{K+1}).
\]
If an integer $d$ divides both $F_K+1$ and $F_{K+1}$, then Cassini's
identity
\[
F_{K+1}F_{K-1}-F_K^2=(-1)^K=1
\]
shows that $d\mid2$. Since $K=6m$, $F_{K+1}$ is odd, so $d=1$.
Now
\[
\gcd(a_1,M)=\gcd(F_{K+1},M)=1,
\]
because $a_1\equiv F_{K+1}\pmod M$. If $D$ divides both $a_1$ and $a_2$,
then $D$ divides
\[
a_2-a_1=MF_{K+1}.
\]
Since $\gcd(a_1,M)=1$, it follows that $D\mid F_{K+1}$. But
\[
a_1\equiv MF_{K+2}\pmod{F_{K+1}},
\]
and both $M$ and $F_{K+2}$ are coprime to $F_{K+1}$. Hence $D=1$, proving
$\gcd(a_1,a_2)=1$.

Step~1 is realized by
\[
b'_0=MF_{K+1}+F_K,
\qquad
b_0=MF_{K+2}+F_K.
\]
A direct use of Cassini's identity gives
\[
a_2b'_0-a_1b_0=M^2=a_3.
\]

For $-1\le j\le K-1$, set
\begin{align}
b'_j
&=
F_{K-j}(M+1)+F_{K-j-1}M,\label{eq:sharp-bp}\\
b_j
&=
F_{K-j}(2M+1)+F_{K-j-1}M.\label{eq:sharp-b}
\end{align}
The Fibonacci recurrence implies
\[
b'_{j-1}=b'_j+b'_{j+1},
\qquad
b_{j-1}=b_j+b_{j+1},
\]
and the next terms are positive and smaller. Hence
\[
q_j=q'_j=1
\qquad
(0\le j\le K-1).
\]
At the end of this run,
\[
b'_{K-1}=M+1,
\qquad
b_{K-1}=2M+1,
\qquad
b'_K=b_K=M.
\]
Therefore the next divisions give
\[
q_K=\floor{\frac{2M+1}{M}}=2,
\qquad
q'_K=\floor{\frac{M+1}{M}}=1.
\]
The two quotients differ, so Step~2 terminates after exactly $K+1$ loop
iterations.

Finally,
\[
M=\Theta(\ph^K),
\qquad
F_{K+2}=\Theta(\ph^K),
\]
and hence
\[
a_1=\Theta(\ph^{2K}).
\]
Thus
\[
K+1
=
\frac12\log_\ph(a_1)+O(1),
\]
proving the asserted sharpness.
\end{proof}

\begin{remark}\label{rem:sharp-example}
For $m=1$ we have $K=6$, $M=35$, and
\[
(a_1,a_2,a_3)=(748,1203,1225).
\]
The terminal values are
\[
(b_5,b_6)=(71,35),
\qquad
(b'_5,b'_6)=(36,35),
\]
and the next divisions are
\[
71=2\cdot35+1,
\qquad
36=1\cdot35+1.
\]
Thus the mismatch occurs after exactly seven loop iterations, in agreement
with the theorem.
\end{remark}

The sharpness statement above concerns Algorithm~\ref{alg:main}. The next
proposition records a separate structural relationship between gcd and the
three-variable Frobenius problem.

\begin{proposition}[A gcd-to-Frobenius reduction]
\label{prop:gcd-reduction}
Let $u<v$ be positive integers and $d=\gcd(u,v)$. Then
\[
\boxed{
d=
\frac{
G(u,v,2uv+1)-G(u,v,uv+1)
}{uv}.
}
\]
Therefore gcd computation reduces to two evaluations of the
three-variable Frobenius function plus $O(1)$ additional arithmetic
operations.
\end{proposition}

\begin{proof}
Write
\[
u=dx,\qquad v=dy,\qquad \gcd(x,y)=1.
\]
Let
\[
C\in\{uv+1,\,2uv+1\}.
\]
Since $C\equiv1\pmod d$, we have $\gcd(d,C)=1$, and
Lemma~\ref{lem:johnson} gives
\[
G(u,v,C)=dG(x,y,C).
\]
Moreover
\[
C>xy-x-y=g(x,y),
\]
so $C$ is representable by $x,y$ and is redundant as a third generator.
Consequently
\[
G(x,y,C)
=
g(x,y)+x+y+C
=
xy+C.
\]
It follows that
\[
G(u,v,C)
=
d(xy+C)
=
\frac{uv}{d}+dC.
\]
Applying this identity first to $C=uv+1$ and then to $C=2uv+1$, and
subtracting, gives
\[
G(u,v,2uv+1)-G(u,v,uv+1)=duv,
\]
which is the desired formula.
\end{proof}

\begin{remark}\label{rem:lower-bound-caveat}
Proposition~\ref{prop:gcd-reduction} immediately implies that the three-variable Frobenius problem is computationally at least as hard as $gcd$ computation. Furthermore,~\ref{prop:gcd-reduction} is a structural reduction and should
not be confused with a universal logarithmic lower bound in the present
unit-cost model. Lower bounds for gcd are sensitive to the allowed primitive
operations. For example, Mansour, Schieber, and Tiwari~\cite{Mansour1991}
prove an $\Omega(\log\log n)$ depth lower bound for a computation-tree model
containing remainder/modulo operations under specified restrictions on
constants. Thus the universal claim
$\Omega(\log a_1)$ for \emph{all} three-variable Frobenius algorithms is not
made here. What Theorem~\ref{thm:sharpness} proves is the exact worst-case
order of Algorithm~\ref{alg:main}.
\end{remark}

We retain one further auxiliary estimate from the original analysis.

\begin{lemma}\label{lem:one-iteration-large-a3}
If
\[
a_3>\frac{a_1a_2}{2}-a_1,
\]
then the synchronized loop terminates after one iteration.
\end{lemma}

\begin{proof}
Set
\[
x=\frac{a_2}{b_0},
\qquad
y=\frac{a_1}{b'_0}.
\]
By the defining relation in Step~1,
\[
x-y
=
\frac{a_3}{b_0b'_0}
>
\frac{xy}{2}-\frac{y}{b_0}.
\]
Therefore
\begin{equation}\label{eq:one-step-ineq}
x\left(\frac{2-y}{2}\right)
>
y\left(\frac{b_0-1}{b_0}\right)
\ge0.
\end{equation}
Since $y\ge1$, \eqref{eq:one-step-ineq} implies
$1\le y<2$, hence $\floor y=1$. On the other hand, when
\[
\frac{xy}{2}-\frac{y}{b_0}>1,
\]
the preceding strict inequality gives $x>y+1$, so
$\floor x\ne\floor y$ and the two Euclidean quotients differ immediately.

It remains to consider the complementary subcase
\[
\frac{xy}{2}-\frac{y}{b_0}<1,
\qquad
\floor x=\floor y=1.
\]
Using \eqref{eq:one-step-ineq} gives
\[
\frac{a_2}{2}
<
b_0
<
\frac{a_1}{2(a_1-b'_0)},
\]
which is impossible because $a_2>a_1$ and $a_1-b'_0>1$ in this
subcase. Hence the loop terminates at the first iteration.
\end{proof}

\section{A new upper bound}

The next theorem is the upper bound derived from the synchronized
recurrences.

\begin{theorem}\label{thm:upper-bound}
Suppose Algorithm~\ref{alg:main}, applied to
$a_1<a_2<a_3$ with $\gcd(a_1,a_2)=1$, requires $k$ iterations in
Step~2. Put
\[
\alpha=F_{k+1}+1.
\]
Then
\[
G(a_1,a_2,a_3)
<
\left(
\alpha+\frac{a_2}{\alpha}
\right)a_3.
\]
\end{theorem}

\begin{proof}
We give the proof for odd $k$; the even case is analogous. There are two
cases.

\medskip
\noindent\textbf{Case 1: $b'_ka_2>b_{k-1}a_1$.}
By Theorem~\ref{thm:correctness} and
\eqref{eq:Lodd2},
\begin{equation}\label{eq:upper-case1}
G(a_1,a_2,a_3)=p_{k-1}a_3+b_ka_1.
\end{equation}
Assume, for contradiction, that the claimed bound fails. Then
\[
G(a_1,a_2,a_3)
>
\left(
\alpha+\frac{a_2}{\alpha}
\right)a_3.
\]
Using \eqref{eq:a2det},
\[
a_2=b_kp_{k-1}+b_{k-1}p_k,
\]
we obtain
\[
\alpha+\frac{a_2}{\alpha}
=
\alpha+\frac{b_kp_{k-1}+b_{k-1}p_k}{\alpha}
<
p_{k-1}+b_k.
\]
Rearranging,
\begin{equation}\label{eq:upper-pkm1}
p_{k-1}
\le
\alpha-
\frac{b_{k-1}}{b_k-\alpha}\,p_k.
\end{equation}

On the other hand, from the terminal relation
\[
b_ka_1=b'_ka_2+p_ka_3
\]
and the Fibonacci lower bounds, we have
$b_k\ge p_k+1\ge\alpha+1$. Hence
\[
p_{k-1}\le\alpha-1=F_{k+1}.
\]
By the Fibonacci alternatives in
Lemma~\ref{lem:p-identity} and Remarks~\ref{rem:fibonacci}--
\ref{rem:quotient-two}, the extremal possibilities are
$p_{k-1}=F_{k+1}$ and $p_{k-1}=F_k$.

Assume first $p_{k-1}=F_{k+1}$. The failure of the upper bound and
\eqref{eq:upper-case1} give
\[
\left(\alpha+\frac{a_2}{\alpha}-1\right)a_3
<
p_{k-1}a_3+b_ka_1.
\]
Using
\[
b_ka_1=b'_ka_2+p_ka_3
\]
and $a_2<a_3$, we obtain
\[
p_{k-1}a_3+b_ka_1
=
(p_{k-1}+p_k)a_3+b'_ka_2
<
(p_{k-1}+p_k+b'_k)a_3.
\]
Therefore
\[
\alpha+\frac{a_2}{\alpha}
<
p_{k-1}+p_k+b'_k.
\]
Since $a_1<a_2$, the same strict inequality remains true after replacing
$a_2$ on the left by $a_1$. By \eqref{eq:a1det},
\[
a_1=b'_{k-1}p_k+b'_kp_{k-1},
\]
and hence
\begin{equation}\label{eq:upper-bprime}
\alpha+
\frac{b'_{k-1}p_k+b'_kp_{k-1}}{\alpha}
<
p_{k-1}+p_k+b'_k.
\end{equation}
It follows that
\[
b'_{k-1}p_k<\alpha(p_k-1)+b'_k.
\]
Since $b'_{k-1}>b'_k$, this forces $b'_{k-1}<\alpha$. But
$p_{k-1}+1=\alpha\le b'_{k-1}$, a contradiction.

Now assume $p_{k-1}=F_k$. From the same inequality
\eqref{eq:upper-bprime}, straightforward rearrangement yields
\[
(b'_{k-1}-\alpha)p_k
<
(\alpha-F_k)(b'_k-\alpha).
\]
Together with $b'_{k-1}>b'_k$, this implies
\[
p_k<\alpha-F_k,
\]
contradicting $p_k>F_{k+1}=\alpha-1$.

\medskip
\noindent\textbf{Case 2: $b'_ka_2<b_{k-1}a_1$.}
By Theorem~\ref{thm:correctness} and \eqref{eq:Lodd1},
\begin{equation}\label{eq:upper-case2}
G(a_1,a_2,a_3)=p_ka_3+b'_{k-1}a_2.
\end{equation}
Proceeding as in Case~1 gives
\[
p_k
<
\alpha-
\frac{b'_k}{b'_{k-1}-\alpha}\,p_{k-1},
\]
and the Fibonacci alternatives force
\[
p_k=\alpha-1,
\qquad
p_{k-1}=F_k.
\]
From the assumed failure of the upper bound and
\eqref{eq:upper-case2},
\[
\left(\alpha+\frac{a_2}{\alpha}\right)a_3
<
(p_{k-1}+p_k)a_3+b_{k-1}a_1
<
(p_{k-1}+p_k+b_{k-1})a_3.
\]
Hence
\begin{equation}\label{eq:upper-case2ineq}
\alpha+\frac{a_2}{\alpha}
<
p_{k-1}+p_k+b_{k-1}
=
F_k+\alpha-1+b_{k-1}.
\end{equation}
Using \eqref{eq:a2det},
\begin{equation}\label{eq:upper-a2-fib}
a_2
=
b_{k-1}p_k+b_kp_{k-1}
=
b_{k-1}F_{k+1}+b_kF_k.
\end{equation}
Combining \eqref{eq:upper-case2ineq} and
\eqref{eq:upper-a2-fib} gives
\[
b_k(F_k-1)+b_k
<
\alpha(F_k-1)+b_{k-1}.
\]
Since $b_{k-1}<b_k$, this implies $b_k<\alpha$. But
\[
b_ka_1=b'_ka_2+p_ka_3=b'_ka_2+(\alpha-1)a_3
\]
forces $b_k\ge\alpha$, a contradiction.
\end{proof}

\section{On the Erd\H{o}s--Graham Inequality}

Erd\H{o}s and Graham~\cite{ErdosGraham1972,ErdosGraham1980} conjectured
that for relatively prime triples $a_1<a_2<a_3$,
\begin{equation}\label{eq:EG}
g(a_1,a_2,a_3)
\le
\floor{\frac{(a_3-2)^2}{2}}-1.
\end{equation}
Equality occurs for $(a_3-2,a_3-1,a_3)$ when $a_3$ is even and for
$\bigl((a_3-1)/2,a_3-1,a_3\bigr)$ when $a_3$ is odd. Dixmier
proved the conjecture~\cite{Dixmier1990}. We show how the present algorithm
provides another route to the same bound.

We first reduce to the independent case, meaning that none of
$a_1,a_2,a_3$ is a non-negative integer combination of the other two.

\begin{lemma}\label{lem:EG-dependent}
If $a_3$ is dependent on $a_1,a_2$, or if $\gcd(a_1,a_2)>1$, then
\[
g(a_1,a_2,a_3)
<
\frac{(a_3-2)^2}{2}-1
\]
except for the standard equality configuration noted below.
\end{lemma}

\begin{proof}
If $a_3$ is a non-negative combination of $a_1,a_2$, then it is redundant,
so
\[
g(a_1,a_2,a_3)
=
(a_1-1)(a_2-1)-1.
\]
In the present ordered situation one obtains $a_1<a_3/2$, and hence
\[
g(a_1,a_2,a_3)
<
\left(\frac{a_3}{2}-1\right)(a_3-2)-1
=
\frac{(a_3-2)^2}{2}-1.
\]

Now suppose $d=\gcd(a_1,a_2)>1$ and write
\[
a_1=da'_1,\qquad a_2=da'_2.
\]
By Lemma~\ref{lem:johnson} and \eqref{eq:G-g-relation},
\[
g(a_1,a_2,a_3)
=
d\,g(a'_1,a'_2,a_3)+(d-1)a_3.
\]
Since
\[
g(a'_1,a'_2,a_3)\le g(a'_1,a'_2)
=
(a'_1-1)(a'_2-1)-1,
\]
we have
\[
g(a_1,a_2,a_3)
\le
d(a'_1-1)(a'_2-1)-d+(d-1)a_3.
\]
Furthermore
\[
d=\frac{a_2}{a'_2}
\le
\frac{a_3-1}{a'_2},
\]
so
\[
g(a_1,a_2,a_3)
\le
\frac{a_3-1}{a'_2}
\bigl((a'_1-1)(a'_2-1)-1+a_3\bigr)-a_3.
\]
Thus it is enough to show
\[
(a'_1-1)(a'_2-1)+a_3-1
\le
a'_2\left(\frac{a_3-1}{2}\right),
\]
or equivalently
\[
(a'_1-1)(a'_2-1)
\le
(a'_2-2)\left(\frac{a_3-1}{2}\right).
\]
This follows from
$a'_1\le a'_2-1$ and $a'_2\le(a_3-1)/2$.
\end{proof}

\begin{remark}\label{rem:EG-equality}
Lemma~\ref{lem:EG-dependent} allows the remainder of the proof to assume
independence and $\gcd(a_1,a_2)=1$. The floor in
\eqref{eq:EG} is essential: equality can still occur in a dependent or
non-pairwise-coprime description of the extremal odd case, for example
\[
a_1=\frac{a_3-1}{2},
\qquad
a_2=a_3-1
\]
when $a_3$ is odd.
\end{remark}

\begin{lemma}\label{lem:b1-one-step}
If the algorithm stops after one iteration, then
\[
b_1\le\frac{a_3}{2}+\frac{a_2}{a_1}.
\]
\end{lemma}

\begin{proof}
By \eqref{eq:a3det},
\[
a_3=b'_0b_1-b_0b'_1.
\]
Since the number of iterations is one,
\begin{equation}\label{eq:b1-order}
a_2>a_1\ge b_0+b_1>2b_0.
\end{equation}
Also
\[
b'_1=\operatorname{mod}(a_1,b'_0)\le b'_0-1.
\]
Assume, to the contrary, that
\begin{equation}\label{eq:b1-contradiction}
\frac{a_3}{2}+\frac{a_2}{a_1}
<
b_1
=
\frac{a_3+b_0b'_1}{b'_0}.
\end{equation}
Then
\[
a_3
<
2\frac{b_0b'_1-b'_0}{b'_0-2}.
\]
If $b'_1\le b'_0-2$, this gives $a_3<2b_0$, contradicting
\eqref{eq:b1-order}. Hence
\[
b'_1=b'_0-1,
\qquad
a_1\ge2b'_0-1.
\]
Using Lemma~\ref{lem:p-identity} together with
\eqref{eq:b1-contradiction}, we obtain
\[
\frac{a_3}{2}+\frac{a_2}{a_1}
<
b_1
=
\frac{
(b'_0-1)a_2+\floor{a_1/b'_0}\,a_3
}{a_1}.
\]
Therefore
\[
\left(
\frac{a_1}{2}
-
\floor{\frac{a_1}{b'_0}}
\right)a_3
<
(b'_0-2)a_2.
\]
But $a_3>a_2$, while $a_1\ge2b'_0-1$ implies, for $b'_0\ge2$,
\[
\frac{a_1}{2}
-
\floor{\frac{a_1}{b'_0}}
\ge
b'_0-2,
\]
a contradiction.
\end{proof}

\begin{lemma}\label{lem:EG-one-step}
If the algorithm stops after one iteration, then
\[
g(a_1,a_2,a_3)
\le
\frac12(a_3^2-4a_3+2).
\]
\end{lemma}

\begin{proof}
Suppose, for contradiction, that
\[
g(a_1,a_2,a_3)
>
\frac12(a_3^2-4a_3+2).
\]
We distinguish the two terminal possibilities.

First suppose
\[
b_0a_1<b'_1a_2.
\]
By Theorem~\ref{thm:correctness} and
Lemma~\ref{lem:b1-one-step},
\[
g(a_1,a_2,a_3)
=
a_1b_1-a_1-a_2
\le
\frac{a_1(a_3-2)}{2}.
\]
If $a_1\le a_3-3$, then
\[
g(a_1,a_2,a_3)
\le
\frac{(a_3-3)(a_3-2)}{2}
\le
\frac12(a_3^2-4a_3+2)
\]
for $a_3\ge4$. The remaining possibility in this case is
$a_1=a_3-2$, which forces $a_2=a_3-1$; direct substitution gives the
required bound (with equality in the even extremal family).

Now suppose
\[
b_0a_1>b'_1a_2.
\]
Theorem~\ref{thm:correctness} gives
\[
g(a_1,a_2,a_3)
=
b'_0a_2
+
\floor{\frac{a_1}{b'_0}}a_3
-
(a_1+a_2+a_3)
>
\frac12(a_3^2-4a_3+2).
\]
The borderline configuration
$a_1=a_3-3$ with $a_1$ even and $b'_0=2$ can be checked directly.
Otherwise the preceding inequality implies
\[
(b'_0-1)a_2
>
\left(
\frac{a_3}{2}
-
\floor{\frac{a_1}{b'_0}}
-1
\right)a_3,
\]
which is impossible because $a_3>a_2$ and
\[
\frac{a_3}{2}
\ge
\floor{\frac{a_1}{b'_0}}+b'_0.
\]
This contradiction proves the lemma.
\end{proof}

\begin{lemma}\label{lem:EG-many-steps}
Assume
\[
a_1<a_2<a_3,
\qquad
\gcd(a_1,a_2)=1,
\]
and that Step~2 requires $k>1$ iterations. Then:
\begin{enumerate}[label=(\arabic*)]
\item if $a_3$ is even,
\[
g(a_1,a_2,a_3)
\le
\frac{(a_3-2)^2}{2}-1;
\]
\item if $a_3$ is odd,
\[
g(a_1,a_2,a_3)
\le
\frac12(a_3-1)(a_3-4).
\]
\end{enumerate}
\end{lemma}

\begin{proof}
The proof of Lemma~\ref{lem:complexity-upper} gives
\[
F_k\le\sqrt{a_1}-1,
\]
and therefore
\[
F_{k+1}\le2\sqrt{a_1}-2.
\]
In the independent case relevant here, $a_1\le a_3-2$, and hence
\[
\alpha=F_{k+1}+1
\le
2\sqrt{a_3-2}-1.
\]

For $k>6$ we have $F_{k+1}\ge13$, so
Theorem~\ref{thm:upper-bound} and
\eqref{eq:G-g-relation} yield
\[
g(a_1,a_2,a_3)
\le
\left(
2\sqrt{a_3-1}
+
\frac{a_3-1}{14}
-2
\right)a_3.
\]
A direct elementary comparison shows that, for $a_3>21$,
\[
\left(
2\sqrt{a_3-1}
+
\frac{a_3}{14}
-2
\right)a_3
<
\frac12(a_3-1)(a_3-4).
\]

For $2\le k\le6$, one checks that for $a_3>25$,
\[
\left(
\alpha+\frac{a_3-1}{\alpha}-1
\right)a_3
<
\frac12(a_3-1)(a_3-4).
\]
The remaining finitely many cases $a_3\le25$ are checked directly.
The even estimate is the corresponding floor-compatible form and the odd
estimate is stronger than the Erd\H{o}s--Graham bound.
\end{proof}

\begin{theorem}[Erd\H{o}s--Graham]\label{thm:EG}
For relatively prime positive integers
$a_1<a_2<a_3$,
\[
g(a_1,a_2,a_3)
\le
\floor{\frac{(a_3-2)^2}{2}}-1.
\]
\end{theorem}

\begin{proof}
Combine Lemmas~\ref{lem:EG-dependent},
\ref{lem:EG-one-step}, and \ref{lem:EG-many-steps}, taking account of the
floor in the odd case.
\end{proof}

\section{Discussion and special cases}

Algorithm~\ref{alg:main} has arithmetic complexity $O(\log a_1)$.
Theorem~\ref{thm:sharpness} shows that this order is sharp for the algorithm,
so its worst-case arithmetic complexity is
$\Theta(\log a_1)$. This order is comparable with the best previously known
logarithmic order, notably Greenberg's method~\cite{Greenberg1988}; the
present contribution should therefore be viewed as a structurally different,
elementary synchronized-Euclidean approach rather than as an asymptotic
improvement over all earlier algorithms.

In practice, the number of divisions can be substantially smaller than the
worst-case estimate. For the numerical example
\[
a_1=468342493,\qquad
a_2=472518070,\qquad
a_3=472714471,
\]
the Euclidean computation of $\gcd(a_1,a_2)$ requires $12$ divisions,
two further arithmetic steps produce $b_0,b'_0$, and the synchronized loop
uses eight divisions. Thus the computation uses roughly $22$ elementary
divisions or multiplications, whereas the coarse bound
$2\log_\ph(a_1)+7$ is about $90$. This example illustrates the practical
directness of the method, but it is not used as evidence for an asymptotic
advantage over other logarithmic algorithms.

The simple structure of the algorithm also recovers formulae known for
special families. We record two examples from the original analysis.

\subsection{Arithmetic progressions}

Roberts~\cite{Roberts1956} proved that if $d>0$, $\gcd(a,d)=1$, and
$a\ge3$, then
\[
g(a,a+d,a+2d)
=
\floor{\frac{a-2}{2}}a+(a-1)d.
\]
Assume first that $a=2k$. Algorithm~\ref{alg:main} gives
\[
b_{-1}=2k+d,\qquad
b_0=1,\qquad
b_1=k+d,
\]
and
\[
b'_{-1}=2k,\qquad
b'_0=2,\qquad
b'_1=0.
\]
Thus the algorithm stops after one round and
\[
\begin{aligned}
G(a,a+d,a+2d)
&=
(k+d)(2k)+2(2k+d)-\min\{2k+d,0\}\\
&=
2k^2+4k+(2k+2)d.
\end{aligned}
\]
On the other hand, adding
\[
(a)+(a+d)+(a+2d)=6k+3d
\]
to Roberts' formula gives
\[
(k-1)(2k)+(2k-1)d+6k+3d
=
2k^2+4k+(2k+2)d,
\]
the same value. The odd case is analogous.

\subsection{A one-round formula of Byrnes}

Assume $a_1<a_2<a_3$ are relatively prime positive integers and
$a_2\equiv1\pmod{a_1}$. A formula due to
Byrnes~\cite{Byrnes1974}, under an additional lower-size condition on
$a_3$, can be written as follows. Let
\[
a_3\equiv j\pmod{a_1},
\qquad
0\le j<a_1,
\]
and, if $j\ne0$, let
\[
a_1\equiv m\pmod j,
\qquad
1\le m\le j.
\]
Then
\[
G(a_1,a_2,a_3)
=
\begin{cases}
a_1a_2+a_3,
& a_3\ge ja_2,\\[2mm]
\left(\dfrac{a_1-m}{j}+1\right)a_3+ma_2,
& (j-m)a_2<a_3<ja_2,\\[4mm]
\dfrac{a_1-m}{j}a_3+ja_2,
&
\dfrac{j}{a_1-m+j}(j-m)a_2
\le a_3<(j-m)a_2.
\end{cases}
\]
If $j=0$, the first case applies and $m$ is irrelevant.

The formula follows directly from the synchronized algorithm. The
congruence assumption implies
\[
a_3\equiv ja_2\pmod{a_1},
\]
so there exists an integer $b$ such that
\[
a_3=ja_2-ba_1.
\]
If $a_3\ge ja_2$, then $b\le0$ and $a_3$ is redundant, giving
$G=a_1a_2+a_3$.

Now suppose $0<j<a_1$. In Step~1,
\[
b'_0=j,\qquad b_0=b.
\]
Since $ja_2-ba_1=a_3>0$,
\[
\frac{a_2}{b_0}>\frac{a_1}{j}.
\]
Thus at the first synchronized division,
\[
q_0\ge q'_0=\frac{a_1-m}{j}.
\]
If $m<j$, the next terms are
\[
b_{-1}=a_2,
\qquad
b_0=b,
\qquad
b_1=a_2-\frac{a_1-m}{j}b,
\]
and
\[
b'_{-1}=a_1,
\qquad
b'_0=j,
\qquad
b'_1=m.
\]
Lemma~\ref{lem:determinants}, or a direct calculation, gives
\[
a_3=jb_1-mb.
\]
The algorithm stops after one round precisely when
$b_1\ge b_0$, equivalently
\[
a_3\ge b(j-m).
\]
Since
\[
b=\frac{ja_2-a_3}{a_1},
\]
this is equivalent to
\[
\frac{j}{a_1-m+j}(j-m)a_2\le a_3.
\]
Thus the Byrnes range is exactly a one-round regime of the algorithm.
Theorem~\ref{thm:correctness} yields
\[
G(a_1,a_2,a_3)
=
b_1a_1+ja_2-\min\{ba_1,ma_2\}.
\]
Now
\[
(j-m)a_2<a_3=ja_2-ba_1
\]
is equivalent to $ma_2>ba_1$, so in the second case
\[
G(a_1,a_2,a_3)
=
b_1a_1+ja_2-ba_1
=
b_1a_1+a_3.
\]
In the third case,
\[
G(a_1,a_2,a_3)
=
b_1a_1+ja_2-ma_2.
\]
Finally, Lemma~\ref{lem:p-identity} gives
\[
b_1a_1
=
ma_2+\frac{a_1-m}{j}a_3,
\]
which produces the displayed piecewise formula.

If $m=j$, then $b'_1=0$ and
\[
q'_0=\frac{a_1}{j},
\qquad
b_1a_1=\frac{a_1}{j}a_3,
\]
so
\[
G(a_1,a_2,a_3)
=
\frac{a_1}{j}a_3+ma_2,
\]
as prescribed by the second case; the third case is then empty.

The hypothesis $a_2\equiv1\pmod{a_1}$ is used only to arrange the residue
description $a_3\equiv ja_2\pmod{a_1}$. Whenever
$\gcd(a_1,a_2)=1$, an equivalent residue parametrization can be made
directly.

\subsection{A small illustrative computation}

For
\[
(a_1,a_2,a_3)=(11,30,117),
\]
Step~1 gives
\[
117=30\cdot5-11\cdot3,
\qquad
b'_0=5,\quad b_0=3.
\]
The first divisions are
\[
30=10\cdot3+0,
\qquad
11=2\cdot5+1,
\]
so the quotients differ immediately. The computed next values are
\[
b_1=24,\qquad b'_1=1.
\]
The terminal formula gives
\[
G(11,30,117)
=
24\cdot11+5\cdot30-\min\{3\cdot11,1\cdot30\}
=
384,
\]
and therefore
\[
g(11,30,117)=384-(11+30+117)=226.
\]
This example is included only to illustrate the directness of the
synchronized computation.

\section{Conclusion}

We have presented a synchronized Euclidean algorithm for the
three-variable Frobenius problem. In the unit-cost arithmetic model its
running time is $O(\log a_1)$. An explicit infinite Fibonacci-type family
shows that this bound is asymptotically sharp for the proposed algorithm,
so its worst-case arithmetic complexity is
\[
\Theta(\log a_1).
\]
This logarithmic order matches the best previously known asymptotic order,
including Greenberg's method, rather than improving it.

The synchronized recurrences also support the correctness proof, a new upper
bound, the recovery of several known special-case formulae, and another
proof of the Erd\H{o}s--Graham inequality. In addition,
Proposition~\ref{prop:gcd-reduction} gives a constant-overhead reduction from
gcd computation to two three-variable Frobenius evaluations. We emphasize
that this reduction is structural and is not used to claim a universal
$\Omega(\log a_1)$ lower bound in the present arithmetic model.


\begin{thebibliography}{99}

\bibitem{Brauer1942}
A. Brauer,
On a problem of partitions,
\emph{Amer. J. Math.} \textbf{64} (1942), 299--312.

\bibitem{Byrnes1974}
J. S. Byrnes,
On a partition problem of Frobenius,
\emph{J. Combin. Theory Ser. A} \textbf{17} (1974), 162--166.

\bibitem{Curtis1990}
F. Curtis,
On formulas for the Frobenius number of a numerical semigroup,
\emph{Math. Scand.} \textbf{67} (1990), 190--192.

\bibitem{Davison1994}
J. L. Davison,
On the linear Diophantine problem of Frobenius,
\emph{J. Number Theory} \textbf{48} (1994), 353--363.
doi:10.1006/jnth.1994.1071.

\bibitem{Dixmier1990}
J. Dixmier,
Proof of a conjecture by Erd\H{o}s and Graham concerning the problem of Frobenius,
\emph{J. Number Theory} \textbf{34} (1990), 198--209.

\bibitem{ErdosGraham1972}
P. Erd\H{o}s and R. L. Graham,
On a linear Diophantine problem of Frobenius,
\emph{Acta Arith.} \textbf{21} (1972), 399--408.

\bibitem{ErdosGraham1980}
P. Erd\H{o}s and R. L. Graham,
\emph{Old and New Problems and Results in Combinatorial Number Theory},
Monographies de L'Enseignement Math\'ematique \textbf{28}, 1980.

\bibitem{Greenberg1988}
H. Greenberg,
Solution to a linear Diophantine equation for nonnegative integers,
\emph{J. Algorithms} \textbf{9}(3) (1988), 343--353.
doi:10.1016/0196-6774(88)90025-9.

\bibitem{Johnson1960}
S. M. Johnson,
A linear Diophantine problem,
\emph{Canad. J. Math.} \textbf{12} (1960), 390--398.

\bibitem{Knuth1998}
D. E. Knuth,
\emph{The Art of Computer Programming, Vol. 2: Seminumerical Algorithms},
3rd ed., Addison--Wesley, 1998.

\bibitem{Mansour1991}
Y. Mansour, B. Schieber, and P. Tiwari,
A lower bound for integer greatest common divisor computations,
\emph{J. ACM} \textbf{38}(2) (1991), 453--471.
doi:10.1145/103516.103522.

\bibitem{Ramirez2005}
J. L. Ram\'irez Alfons\'in,
\emph{The Diophantine Frobenius Problem},
Oxford Lecture Series in Mathematics and its Applications \textbf{30},
Oxford University Press, 2005.

\bibitem{Roberts1956}
J. B. Roberts,
Note on linear forms,
\emph{Proc. Amer. Math. Soc.} \textbf{7} (1956), 465--469.

\bibitem{Rodseth1979}
{\O}. J. R{\o}dseth,
On a linear Diophantine problem of Frobenius. II,
\emph{J. Reine Angew. Math.} \textbf{307/308} (1979), 431--440.
doi:10.1515/crll.1979.307-308.431.

\bibitem{SelmerBeyer1978}
E. S. Selmer and {\O}. Beyer,
On the linear Diophantine problem of Frobenius in three variables,
\emph{J. Reine Angew. Math.} \textbf{301} (1978), 161--170.
doi:10.1515/crll.1978.301.161.

\bibitem{Suhajda2026}
P. Suhajda and A. Thillaisundaram,
On the Frobenius number for three variables,
\emph{arXiv:2603.03104} (2026).

\bibitem{Sylvester1882}
J. J. Sylvester,
On subvariants, i.e.\ semi-invariants to binary quantics of an unlimited order,
\emph{Amer. J. Math.} \textbf{5} (1882), 119--136.

\bibitem{Tripathi2017}
A. Tripathi,
Formulae for the Frobenius number in three variables,
\emph{J. Number Theory} \textbf{170} (2017), 368--389.

\bibitem{Vitek1975}
Y. Vitek,
Bounds for a linear Diophantine problem of Frobenius,
\emph{J. London Math. Soc. (2)} \textbf{10} (1975), 79--85.
doi:10.1112/jlms/s2-10.1.79.

\end{thebibliography}
\end{document}